\documentclass[12pt]{article}
\usepackage[centertags]{amsmath}
\usepackage{amsfonts}
\usepackage{amssymb}
\usepackage{amsthm}
\usepackage{newlfont}
\usepackage{graphicx}

\newfont{\bb}{msbm10 at 12pt}
\def\r{\hbox{\bb R}}

\def\e{\hbox{\bf E}}
\def\t{\hbox{\bf T}}
\def\n{\hbox{\bf N}}
\def\b{\hbox{\bf B}}

\newtheorem{theorem}{Theorem}[section]

\begin{document}

\title{Determination of the position vectors of general helices from intrinsic equations in $\e^3$ }
\author{ Ahmad T. Ali\\Mathematics Department\\
 Faculty of Science, Al-Azhar University\\
 Nasr City, 11448, Cairo, Egypt\\
email: atali71@yahoo.com}

\maketitle
\begin{abstract} In this paper, we prove that the position vector of every space curve satisfies a vector differential equation of fourth order. Also, we determine the parametric representation of the position vector $\psi=\Big(\psi_1,\psi_2,\psi_3\Big)$ of general helices from the intrinsic equations $\kappa=\kappa(s)$ and $\tau=\tau(s)$ where $\kappa$ and $\tau$ are the curvature and torsion of the space curve $\psi$, respectively. Our result extends some knwown results. Moreover, we give four examples to illustrate how to find the position vector from the intrinsic equations of general helices.

\end{abstract}



\section{Introduction }
{\it Helix} is one of the most fascinating curves in science and nature. Scientist have long held a fascinating, sometimes bordering on mystical obsession, for helical structures in nature. Helices arise in nano-springs, carbon nano-tubes, $\alpha$-helices, DNA double and collagen triple helix, lipid bilayers, bacterial flagella in salmonella and escherichia coli, aerial hyphae in actinomycetes, bacterial shape in spirochetes, horns, tendrils, vines, screws, springs, helical staircases and sea shells \cite{choua, lucas, watson}. Also we can see the helix curve or helical structures in fractal geometry, for instance hyperhelices \cite{toledo}. In the field of computer aided design and computer graphics, helices can be used for the tool path description, the simulation of kinematic motion or the design of highways, etc. \cite{yang}. From the view of differential geometry, a helix is a geometric curve with non-vanishing constant curvature $\kappa$ and non-vanishing constant torsion $\tau$ \cite{barros}. The helix may be called a {\it circular helix} or {\it W-curve} \cite{ilarslan, mont1}.

Indeed a helix is a special case of the {\it general helix}. A curve of constant slope or general helix in Euclidean 3-space $\e^3$ is defined by the property that the tangent makes a constant angle with a fixed straight line called the axis of the general helix. A classical result stated by Lancret in 1802 and first proved by de Saint Venant in 1845 (see \cite{struik} for details) says that: {\it A necessary and sufficient condition that a curve be a general helix is that the ratio $$\dfrac{\kappa}{\tau}$$ is constant along the curve, where $\kappa$ and $\tau$ denote the curvature and the torsion, respectively}.

A general helices or {\it inclined curves} are well known curves in classical differential geometry of space curves \cite{milm} and we refer to the reader for recent
works on this type of curves \cite{barros, ba1, gl1, mont2, sc, tur2}. Many important results in the theory of the curves in $\e^3$ were initiated by G. Monge and G. Darboux pioneered the moving frame idea. Thereafter, F. Frenet defined his moving frame and his special equations which play important role in mechanics and kinematics as well as in differential geometry \cite{boyer}. For unit speed curve with non-vanishing curvature $\kappa \neq 0$, it is well-known the following result \cite{hacis}:

\begin{theorem}\label{th-main} A curve is defined uniquely by its curvature and torsion as function of a natural parameters.
\end{theorem}

The equations
$$
\kappa=\kappa(s),\,\,\,\,\,\tau=\tau(s)
$$
which give the curvature and torsion of a curve as functions of $s$ are called the {\it natural} or {\it intrinsic equations} of a curve, for they completely define the curve.

Given two functions of one parameter (potentially curvature and torsion parameterized by arc-length) one might like to find an arc-length parameterized curve for which the two functions work as the curvature and the torsion. This problem, known as {\it solving natural equations}, is generally achieved by solving a {\it Riccati equation} \cite{struik}. Barros et. al. \cite{ba2} showed that the general helices in Euclidean 3-space $\e^3$ and in the three-sphere $\bold{S}^3$ are geodesic either of right cylinders or of Hopf cylinders according to whether the curve lies in $\e^3$ or $\bold{S}^3$, respectively.

In classical differential geometry, The problem of the determination of parametric representation of the position vector of an arbitrary space curve according to the intrinsic equations is still open \cite{eisenh, lips}. This problem is solved in the case of a plane curve $(\tau=0)$ and in the case of circular helix ($\kappa$ and $\tau$ are both non-vanishing constants). However, This problem is not solved in the case of the general helix ($\dfrac{\tau}{\kappa}$ is constant).

Our main result in this work is to proven that the components of the position vector of every space curve satisfies a vector differential equation of forth order and determined the parametric representation of the  position vector $\psi$ from intrinsic equations in $\e^3$ for a general helix $\dfrac{\tau}{\kappa}=\cot[\alpha]$, where the constant $\alpha$ is the angle between the tangent of the curve $\psi$ and the constant vector $\bold{U}$ called the axis of a general helix.

\section{Preliminaries }

In Euclidean space $\e^3$, it is well known that to each unit speed curve with at least four continuous derivatives, one can associate three mutually orthogonal unit vector fields $\t$, $\n$ and $\b$ are respectively, the tangent, the principal normal and the binormal vector fields \cite{hacis}.

We consider the usual metric in Euclidean 3-space $\e^3$, that is,
$$
\langle,\rangle=dx_1^2+dx_2^2+dx_3^2,
$$
where $(x_1,x_2,x_3)$ is a rectangular coordinate system of $\e^3$.  Let $\psi:I\subset\r\rightarrow\e^3$, $\psi=\psi(s)$, be an arbitrary curve in $\e^3$. The curve $\psi$ is said to be of unit speed (or parameterized by the  arc-length) if $\langle\psi'(s),\psi'(s)\rangle=1$ for any $s\in I$. In particular, if $\psi(s)\not=0$ for any $s$, then it is possible to re-parameterize $\psi$, that is, $\alpha=\psi(\phi(s))$ so that $\alpha$ is parameterized by the arc-length. Thus, we will assume throughout this work that $\psi$ is a unit speed curve.

Let $\{\t(s),\n(s),\b(s)\}$ be the moving frame along $\psi$, where the vectors $\t, \n$ and $\b$ are mutually orthogonal vectors satisfying $\langle\t,\t\rangle=\langle\n,\n\rangle=\langle\b,\b\rangle=1$.
The Frenet equations for $\psi$ are given by (\cite{struik})
\begin{equation}\label{u2}
 \left[
   \begin{array}{c}
     \t' \\
     \n' \\
     \b' \\
   \end{array}
 \right]=\left[
           \begin{array}{ccc}
             0 & \kappa & 0 \\
             -\kappa & 0 & \tau \\
             0 & -\tau & 0 \\
           \end{array}
         \right]\left[
   \begin{array}{c}
     \t \\
     \n \\
     \b \\
   \end{array}
 \right].
 \end{equation}

If $\tau(s)=0$ for any $s\in I$, then $\b(s)$ is a constant vector $V$ and the curve $\psi$ lies in a $2$-dimensional affine subspace orthogonal to $V$, which is isometric to the Euclidean $2$-space $\e^{2}$.

We observe that the Frenet equations form a system of three vector differential equations of the first order in $\t, \n$ and $\b$. It is reasonable to ask, therefore, given arbitrary continuous functions $\kappa$ and $\tau$, whether or not there exist solutions $\t, \n, \b$ of the Frenet equations, and hence, since $\psi^{\prime}=\t$, a curve
$$
\psi=\int\,\t\,ds+\bold{C}
$$
which the prescribed curvature and torsion. The answer is in the affirmative and is given by

\begin{theorem}{\bf (Fundamental existence and uniqueness theorem for space curve).}\label{th-main} Let $\kappa(s)$ and  $\tau(s)$ be arbitrary continuous function on $a\leq s \leq b$. Then there exists, except for position in space, one and only one curve $C$ for which $\kappa(s)$ is the curvature, $\tau(s)$ is the torsion and $s$ is a natural parameter along $C$.
\end{theorem}

\section{Position vectors of space curves }

\begin{theorem}\label{th-main} Let $\psi=\psi(s)$ be an unit speed curve. Then, position $\psi$ satisfies a vector differential forth order as follows
\begin{equation}\label{u21}
\dfrac{d}{ds}\Big[\dfrac{1}{\tau}\dfrac{d}{ds}\Big(\dfrac{1}{\kappa}\dfrac{d^2\psi}{ds^2}\Big)\Big]+
\Big(\dfrac{\kappa}{\tau}+\dfrac{\tau}{\kappa}\Big)\dfrac{d^2\psi}{ds^2}+
\dfrac{d}{ds}\Big(\dfrac{\kappa}{\tau}\Big)\dfrac{d\psi}{ds}=0.
\end{equation}
\end{theorem}

{\bf Proof.} Let $\psi=\psi(s)$ be an unit speed curve. If we substitute the first equation of (\ref{u2})
to the second equation of (\ref{u2}), we have
\begin{equation}\label{u3}
\b=\dfrac{1}{\tau}\dfrac{d}{ds}\Big(\dfrac{1}{\kappa}\dfrac{d\t}{ds}\Big)+\dfrac{\kappa}{\tau}\t.
\end{equation}
The last equation of (\ref{u2}) takes the form
\begin{equation}\label{u4}
\dfrac{d}{ds}\Big[\dfrac{1}{\tau}\dfrac{d}{ds}\Big(\dfrac{1}{\kappa}\dfrac{d\t}{ds}\Big)\Big]+
\Big(\dfrac{\kappa}{\tau}+\dfrac{\tau}{\kappa}\Big)\dfrac{d\t}{ds}+
\dfrac{d}{ds}\Big(\dfrac{\kappa}{\tau}\Big)\t=0.
\end{equation}
Denoting $\dfrac{d\psi}{ds}=\t$, we have a vector differential equation of fourth order (\ref{u21}) as desired.

The equation (\ref{u4}) can be written in the following simple form:
\begin{equation}\label{u5}
\dfrac{d}{d\theta}\Big(f\,\dfrac{d^2\t}{d\theta^2}\Big)+
\Big(\dfrac{f^2+1}{f}\Big)\dfrac{d\t}{d\theta}+
\dfrac{df}{d\theta}\t=0,
\end{equation}
where $f=f(\theta)=\dfrac{\kappa(\theta)}{\tau(\theta)}$ and $\theta=\int\kappa(s)ds$. By means of solution of the above equation, position vector of an arbitrary space curve can be determined. However, for general helices, we have

\begin{theorem}\label{th-main2} The position vector of a general helix are computed in the natural parameter form
\begin{equation}\label{u211}
\psi(s)=\sin[\alpha]\int\Big(\cos\Big[\csc[\alpha]\int\kappa(s)ds\Big],\sin\Big[\csc[\alpha]\int\kappa(s)ds\Big], \cot[\alpha]\Big)ds+\bold{C}
\end{equation}
or in the parametric form
\begin{equation}\label{u212}
\psi(\phi)=\int\dfrac{\sin^2[\alpha]}{\kappa(\phi)}\Big(\cos[\phi],\sin[\phi], \cot[\alpha]\Big)d\phi+\bold{C},\,\,\,\phi=\csc[\alpha]\int\kappa(s)ds.
\end{equation}
\end{theorem}

{\bf Proof:} If $\psi$ is a general helix whose tangent vector $\psi^{\prime}$ makes an angle $\alpha$ with the axis $U$, then we can write $f(\theta)=\tan[\alpha]$. Therefore the equation (\ref{u212}) becomes
\begin{equation}\label{u6}
\dfrac{d^3\t}{d\theta^3}+\csc^2[\alpha]\dfrac{d\t}{d\theta}=0.
\end{equation}
or
\begin{equation}\label{u61}
\dfrac{d^3\t}{d\phi^3}+\dfrac{d\t}{d\phi}=0, \,\,\,\,\,\phi=\csc[\alpha]\theta.
\end{equation}

If we write the tangent vector $\t=\Big(T_1, T_2, T_3\Big)$ the general solution of (\ref{u61}) takes the form
\begin{equation}\label{u7}
\t(\phi)=T_i(\phi)\bold{e}_i=\Big(a_i\cos[\phi]+b_i\sin[\phi]+
c_i\Big)\bold{e}_i,\,\,\,i=1,2,3,
\end{equation}
where $a_i, b_i, c_i\in R$ for i=1,2,3.

Hence the curve $\psi$ is general helix, i.e. the tangent vector $\t$ makes an constant angle $\alpha$ with the constant vector called the axis of the helix. So, with out loss of generality, we take the axis of helix is parallel to $\bold{e}_3$. Then $T_3=\langle\t,\bold{e}_3\rangle=\cos[\alpha]$ which leads to $a_3=b_3=0$ and $c_3=\cos[\alpha]$.

On other hand the tangent vector $\t$ is a unit vector, so the following condition is satisfied
\begin{equation}\label{u8}
T_1^2+T_2^2+T_3^2=1,
\end{equation}
which leads to
\begin{equation}\label{u9}
\begin{array}{ll}
&\Big(a_1\cos[\phi]+b_1\sin[\phi]+c_1\Big)^2+\Big(a_2\cos[\phi]+b_2\sin[\phi]+
c_2\Big)^2=\sin^2[\alpha].
\end{array}
\end{equation}
The above equation can be written in the form
\begin{equation}\label{u10}
A_0+\sum_{i=1}^2\Big[A_i\cos[i\,\phi]+B_i\sin[i\,\phi]\Big]=0,
\end{equation}
where
\begin{equation}\label{u11}
\left\{\begin{array}{ll}
A_2&=\dfrac{1}{2}\Big(a_1^2+b_1^2-a_2^2-b_2^2\Big)\\
B_2&=a_1a_2+b_1b_2\\
A_1&=2(a_1a_3+b_1b_3)\\
B_1&=2(a_2a_3+b_2b_3)\\
A_0&=\dfrac{1}{2}\Big[a_1^2+b_1^2+a_2^2+b_2^2+2\Big(a_3^2+b_3^2-\sin^2[\alpha]\Big)\Big].
\end{array}\right.
\end{equation}
If equation (\ref{u10}) is satisfied, then all coefficients must be zero, so we have the following set of algebraic equations in the six unknowns $a_1, a_2, a_3, b_1, b_2$ and $b_3$.
\begin{equation}\label{u12}
A_i=0,\,\,\, \forall\,\,\, i=1,2,...,5.
\end{equation}
Solving the five algebraic equations above we obtain four cases of solutions as the following:
\begin{equation}\label{u13}
\left\{\begin{array}{ll}
a_3=b_3=0,\,\,\,\,b_2=a_1,\,\,\,\, a_2=-b_1,\,\,\,\,b_1=\pm\sqrt{\sin^2[\alpha]-a_1^2}\\
a_3=b_3=0,\,\,\,\,b_2=-a_1,\,\,\,\, a_2=b_1,\,\,\,\,b_1=\pm\sqrt{\sin^2[\alpha]-a_1^2}.
\end{array}\right.
\end{equation}
The four cases above leas to the one general form solution, so the equation (\ref{u7}) takes the form:
\begin{equation}\label{u14}
\begin{array}{ll}
\t(\phi)=&\Big(a_1\cos[\phi]-\sqrt{\sin^2[\alpha]-a_1^2}\sin[\phi]\Big)\bold{e}_1\\
&+
\Big(\sqrt{\sin^2[\alpha]-a_1^2}\cos[\phi]+a_1\sin[\phi]\Big)\bold{e}_2+
\cos[\alpha]\bold{e}_3,
\end{array}
\end{equation}
or in the following form:
\begin{equation}\label{u15}
\begin{array}{ll}
\t(\phi)=\Big(\sin[\alpha]\cos[\phi+\varepsilon],
\sin[\alpha]\sin[\phi+\varepsilon], \cos[\alpha]\Big).
\end{array}
\end{equation}
where $\varepsilon=\arctan\Big[\sqrt{\dfrac{\sin^2[\alpha]}{a_1^2}-1}\Big]$. Without loss of generality we can written:
\begin{equation}\label{u16}
\begin{array}{ll}
\t(\phi)=\Big(\sin[\alpha]\cos[\phi],
\sin[\alpha]\sin[\phi], \cos[\alpha]\Big).
\end{array}
\end{equation}
By Integrating  the above equation with respect to $s$ along with $\phi=\csc[\alpha]\int\kappa(s)ds$, we have the two equations (\ref{u211}) and (\ref{u212}) which it completes the proof.

\section{Examples}

In this section, we take several choices for the curvature $\kappa$ and torsion $\tau$, and next, we apply Theorem \ref{th-main2}.

{\bf Example 1.} The case of plane curve $\tau=0, \kappa=\kappa(s)$, i.e., $\alpha=\dfrac{\pi}{2}$. Then the tangent vector takes the form:
\begin{equation}\label{u17}
\begin{array}{ll}
\t(\phi)=\Big(\cos[\phi],
\sin[\phi], 0\Big).
\end{array}
\end{equation}
Integrate the above equation with respect to $s$ along with $\phi=\int\kappa(s)ds$, we have the natural form of the plane curve as the following:
\begin{equation}\label{u213}
\psi(s)=\int\Big(\cos\Big[\int\kappa(s)ds\Big], \sin\Big[\int\kappa(s)ds\Big], 0\Big)ds+\bold{C}
\end{equation}
or in the parametric form
\begin{equation}\label{u214}
\psi(\phi)=\int\dfrac{1}{\kappa(\phi)}\Big(\cos[\phi], \sin[\phi], 0\Big)d\phi+\bold{C},\,\,\,\phi=\int\kappa(s)ds,
\end{equation}
which is the well-known equation of a plane curve with an arbitrary curvature $\kappa=\kappa(s)$.

{\bf Example 2.} The case of a curve when both of the curvature and torsion are constants, i.e., $\kappa=\dfrac{\sin[\alpha]}{a},\,\tau=\dfrac{\cos[\alpha]}{a}$. Then position vector takes the form:
\begin{equation}\label{u18}
\begin{array}{ll}
\psi=\sin[\alpha]\int\Big(\cos\Big[\dfrac{s}{a}\Big],
\sin\Big[\dfrac{s}{a}\Big],\cot[\alpha]\Big)ds+\bold{C}.
\end{array}
\end{equation}
Integrating the above equation and putting $s=a\,\phi$, we obtain the parametric representation of this curve as the following:
\begin{equation}\label{u191}
\begin{array}{ll}
\psi=a\,\sin[\alpha]\Big(\sin[\phi],-\cos[\phi], \cot[\alpha]\,\phi\Big)+\bold{C}.
\end{array}
\end{equation}
which is the equation of the circular helix.

{\bf Example 3.} The case of a general helix with $\kappa=\dfrac{\sin[\alpha]}{a\,s}$ and $\tau=\dfrac{\cos[\alpha]}{a\,s}$. Then position vector takes the form:
\begin{equation}\label{u20}
\begin{array}{ll}
\psi=a\sin[\alpha]\int\exp[a\phi]\Big(\cos[\phi],
\sin[\phi], \cot[\alpha]\Big)d\phi+\bold{C}.
\end{array}
\end{equation}
Integrating the above equation, we obtain the parametric representation of this curve as the following:
\begin{equation}\label{u21}
\begin{array}{ll}
\psi=\dfrac{a\sin[\alpha]}{1+a^2}\exp[a\,\phi]\Big(
\sin[\phi]+a\cos[\phi], a\sin[\phi]-\cos[\phi], \dfrac{(1+a^2)\cot[\alpha]}{a}\Big)+\bold{C},
\end{array}
\end{equation}
which is the equation of a helix on a cone of revolution.

{\bf Example 4.} The case of a general helix with $\kappa=\dfrac{a\,\sin[\alpha]}{a^2+s^2}$ and $\tau=\dfrac{a\,\cos[\alpha]}{a^2+s^2}$. Then position vector takes the form:
\begin{equation}\label{u20}
\begin{array}{ll}
\psi=a\sin[\alpha]\int\sec^2[\phi]\Big(\cos[\phi],
\sin[\phi], \cot[\alpha]\Big)d\phi+\bold{C}.
\end{array}
\end{equation}
After some computations, we can obtain the parametric representation of this curve as the following:
\begin{equation}\label{u21}
\begin{array}{ll}
\psi=a\,\sin[\alpha]\Big(
\theta, \cosh[\theta], \cot[\alpha]\sinh[\theta]\Big)+\bold{C},
\end{array}
\end{equation}
where $\theta=\sinh^{-1}[\tan[\phi]]$.


\end{document}